\documentclass[11pt]{article}

\usepackage[latin1]{inputenc}
\usepackage{graphicx}
\usepackage{epsfig}
\newenvironment{proof}{\trivlist
\item[\hskip\labelsep{\it Proof}\,:]}{\hfill{$q.e.d$}\endtrivlist}

\newtheorem{theorem}{Theorem}[section]

\newtheorem{proposition}{Proposition}[section]
\newtheorem{corollary}{Corollary}[section]

\newtheorem{remark}{Remark}[section]

\newfont{\bb}{msbm10 at 10pt}
\def\r{\hbox{\bb R}}
\def\l{\hbox{\bb L}}
\def\e{\hbox{\bf e}}

\begin{document}

\title{Surfaces of annulus type with constant mean curvature\\
in Lorentz-Minkowski space}
\author{Rafael L\'opez\\
Departmento de Geometría y Topología\\
Universidad de Granada\\
18071 Granada (Spain)\\
e-mail:{\tt rcamino@ugr.es}}

\date{}

\maketitle

\begin{abstract}
In this paper we solve the Plateau problem for spacelike surfaces with constant mean curvature in
Lorentz-Minkowski three-space $\l^3$ and  spanning two circular (axially symmetric) contours in parallel planes. 
We prove that rotational symmetric surfaces are the only  compact spacelike   surfaces in $\l^3$ 
of  constant mean curvature bounded by 
two concentric circles in parallel planes. As conclusion, we characterize spacelike surfaces of revolution with constant mean curvature 
as the only that either i) are the solutions of the exterior Dirichlet problem for constant boundary data or ii) have an isolated 
conical-type singularity.
%\subclass{53C42, 53B30, 53C50}
\end{abstract}

%%%%%%%%%%%%%%%%%%%%%%%%%%%%%%
\section{Introduction and statement of the results}
%%%%%%%%%%%%%%%%%%%%%%%%%%%%%%

Let $\l^3$ denote the 3-dimensional Lorentz-Minkowski space, that is, the real vector space $\r^3$ endowed with the 
Lorentzian metric $\langle,\rangle=dx_1^2+dx_2^2-dx_3^2$, where $x=(x_1,x_2,x_3)$ are the canonical 
coordinates in $\l^3$. An  immersion $x: \Sigma\rightarrow \l^3$ of a smooth surface 
$\Sigma$  is called spacelike if the induced metric on the surface is positive definite. In this setting, the notions of 
the first and  second fundamental form, and the mean curvature are defined in the same way as 
 in Euclidean space. This article  deals with spacelike immersed surfaces $x:\Sigma 
\rightarrow\l^3$ with {\it constant mean curvature} $H$. From the variational viewpoint, it is well known that such surfaces are critical points of the area functional for variations which preserve a suitable volume function. 

Constant mean curvature spacelike submanifolds of a Lorentzian manifold have interest in relativity theory. In this 
setting, there is interest of finding 
 real-valued functions on a given spacetime, all of whose level sets have constant mean curvature. The mean curvature function may then be used as a global time coordinate and provide a time gauge which  have already been many applications. For example, they have been used to prove positivity of mass \cite{sy}, analyse the space of solutions of Einstein equations \cite{fmm} and in numerical integration schemes for Einstein equations \cite{es,pi}. Further references can be found in review papers such as \cite{cy,mt}.  

A natural family to study consists of the surfaces of revolution. By a such surface, we
mean a surface that it is invariant by the action of a uniparametric subgroup of 
isometries of $\l^3$. The isometries group of $\l^3$ is the semi-direct product of the 
translations group and the orthogonal Lorentzian group $O(1,2)$. With respect to the 
orthogonal group, there are three one-parameter subgroups of isometries of $\l^3$ depending on the 
causal character of the axis.  In this paper we are interested by those surfaces whose axis is a timelike line. 
We call also these surfaces
rotational symmetric surfaces. After a 
Lorentzian transformation, we can suppose that the axis is a vertical line and so, the surface
is foliated by Euclidean circles in horizontal planes and centered at the axis. 

Recall that  in the Euclidean 
setting, the  surfaces of revolution with constant mean curvature  were characterized by Delaunay  in 1841 as follows: 
their profile curves are
obtained by rolling a given conic section on a line in a plane and rotating about that line the trace of a focus \cite{de}. 
In Lorentzian-Minkowski space $\l^3$,  spacelike surfaces of revolution with constant mean curvature were also 
characterized  by results of the same kind \cite{hn}.

Rotational symmetric surfaces with constant mean curvature have an 
important role in the study of spacelike constant mean curvature of $\l^3$ since
they can be used as   barrier surfaces. For example,  this occurs in 
the general scheme in the  solvability of the Dirichlet problem for the
mean curvature equation, by establishing the necessary 
$C^0$ and $C^1$ estimates (for a  general guide we refer to  \cite{gt} and \cite{b1,bs,t} in this context). They are also useful in the 
study of the singularities of a (weakly) spacelike surface. Singularities appear by the  degeneracy of the ellipticity of the mean curvature equation that can drop the regularity of the metric. Rotational symmetric surfaces with constant mean curvature 
allow to control the geometry of singularities. For example, Bartnik  \cite{b2} proved that an isolated singular point
in a spacelike surface in $\l^3$ with $C^1$ mean curvature corresponds to a regular point of the surface or a point where the surface is asymptotic to the light cone at this point. In the last case, the point is called a conical-type singularity (see \cite{ec} for maximal surfaces).

We  seek spacelike surfaces with constant mean curvature and spanning two concentric circles lying in parallel planes. Since the boundary of the surface is  rotational symmetric, it is natural to expect the existence of such a surface 
with rotational  symmetry.  However, a standard application of the Alexandrov reflection method cannot prove that such 
surface inherits the symmetries of the boundary. For example, we would need that the surface lies  in the slab determined by 
the two parallel planes (see \cite[Th. 11]{alp}). Indeed, one cannot expect that this occurs because there examples that
show the contrary: see Figure \ref{figure41}.

We may assume without loss of generality that the circles lie in planes parallel to the plane $\{x_3=0\}$. We introduce the following notation. For real numbers $a$ and $r>0$, let 
$$\Gamma(r,a)=\{(r\cos{\theta},r\sin{\theta},a)\in\r^3;0\leq\theta\leq 2\pi\}.$$
We then pose the following

\begin{quote} {\it Problem:} Given real constants $r,R> 0$ and $a, b, H\in\r$, under what condition on $r, R, a, b$ and $H$ there exists
an annulus-type spacelike surface spanning $\Gamma(r,a)\cup\Gamma(R,b)$
  and with constant mean curvature $H$.
\end{quote}
For example, the spacelike property of the surface imposes that $r\not= R$ (we will suppose
that $r<R$).
Bartnik and Simon \cite{bs} have shown that the Dirichlet problem for the constant  mean curvature equation for $H$
can be solved with merely 
the existence of a spacelike surface  spanning the boundary values. In our case, this means that
$|a-b|<R-r$. 

  In this paper we study the existence of a spacelike rotational symmetric  surface spanning $\Gamma(r,a)\cup\Gamma(R,b)$
 of the form
\begin{equation}\label{form}
X(t,\theta)=(t \cos{\theta},t\sin{\theta}, f(t)),\hspace*{1cm}
t\in [r,R],0\leq \theta\leq 2\pi
\end{equation}
where $f\in C^0([r,R])\cap C^{\infty}(r,R)$, and with boundary values
$$f(r)=a,\hspace*{1cm}f(R)=b.$$

Our first result says us that the same condition as  in \cite{bs} assures  the existence of rotational graphs 
bounded by  $\Gamma(r,a)\cup\Gamma(R,b)$.
Exactly, we have

\begin{theorem}\label{t1} Let  $0<r<R<\infty$ and $a, b\in\r$. Then the following conditions are equivalent:
\begin{enumerate}
\item[(i)] There is a rotational spacelike surface of the form  (\ref{form}) with constant mean curvature 
and spanning $\Gamma(r,a)\cup\Gamma(R,b)$
\item[(ii)] The numbers $r, R, a $ and $b$ satisfy the condition
$$\frac{|a-b|}{R-r}<1.$$
\end{enumerate} 
Moreover, we have the following properties:
\begin{enumerate}
\item The surface can extend to be a graph over $\r^2\setminus\{(0,0)\}$. 
\item At the origin, the surface has a singularity of conical-type, except when  the surface is a 
horizontal plane or a hyperbolic plane. 
\item If $H\not=0$, the surface is asymptotic to a light cone at infinity.
\end{enumerate}
\end{theorem}

In particular, the uniqueness of the Dirichlet problem implies that 
Bartnik-Simon solutions are surfaces of revolution provided that the boundary is rotational symmetric. 
Since spacelike compact surfaces are essentially graphs, we conclude (See Corollary \ref{co})  

\begin{quote}{\it Surfaces of revolution are the only compact spacelike  surfaces  in $\l^3$ with  constant mean curvature 
bounded by two concentric circles in parallel planes.}
\end{quote}

We end this article with two results that characterize the surfaces of revolution in the family of spacelike surfaces with constant 
mean curvature. First we show the uniqueness of the Dirichlet problem for the exterior of a disk:

\begin{theorem}\label{t2} Let $u=u(x_1,x_2)$ define a spacelike surface with constant mean 
curvature in the domain $\Omega=\{(x_1,x_2)\in\r^2;x_1^2+x_2^2>r^2\}$, $r>0$, such that 
$u=a\in\r$ on $\partial\Omega$. Then $u$ describes a surface of revolution.
\end{theorem}

The second result concerns with spacelike surfaces having an isolated singularity. 
A result due to   Ecker shows  that
the Lorentzian catenoids of $\l^3$ are the only entire maximal surfaces with an  isolated singularity \cite[Th.  1.6 ]{ec}. When 
the mean curvature is a non-zero constant, we prove

\begin{theorem}\label{t3}
Entire spacelike constant mean curvature surfaces in Minkowski space $\l^3$ having an isolated singularity are, up to Lorentz-transformations, surfaces of revolution.
\end{theorem}

The proof of the last two theorems involves the study of the flux of a closed curve in a 
spacelike surface, together an application of the maximum principle for surfaces with constant 
mean curvature. 

This paper consists of four sections. Section 2 is a preparatory section where we will mention basic properties of the compact spacelike surfaces with constant mean curvature. Section 3 will be devoted to prove Theorem 
\ref{t1} by analysing the different cases that appear.  
Last, Theorems \ref{t2} and \ref{t3} will be proved in Section 4.

%%%%%%%%%%%%%%%%%%%%%%%%%
\section{Geometric preliminaries}
%%%%%%%%%%%%%%%%%%%%%%%%%

Let $x:\Sigma\rightarrow\l^3$ be a smooth spacelike immersion of a surface $\Sigma$. 
Observe that $\e_{\bf 3}=(0,0,1)\in\l^3$ is a unit timelike vector field globally defined on $\l^3$, which determines 
a time-orientation on $\l^3$. This allows us to choose a unique unit normal vector field $N$ on $\Sigma$ which is in the 
same time-orientation as $\e_3$, and hence we may assume that $\Sigma$ is oriented by $N$. We will refer to 
$N$ as the future-directed Gauss map of $\Sigma$. In this 
article  all spacelike surfaces  will be oriented according to this orientation.

In Lorentz-Minkowski space there are not closed spacelike surface. Thus, any 
compact spacelike surface   has non-empty 
boundary. If $\Gamma$ is  a closed curve in $\l^{3}$ and  $x: \Sigma\rightarrow\l^{3}$ is a 
spacelike immersion of a compact surface, we say that the boundary of $\Sigma$ is  $\Gamma$ if the restriction  $x:\partial  \Sigma
\rightarrow \Gamma$ is a diffeomorphism. 
For spacelike surfaces, the projection $\pi:\l^3\rightarrow \Pi=\{x_3=0\}$, $\pi(x_1,x_2,x_3)=(x_1,x_2,0)$ is  a local diffeomorphism 
between $\mbox{int }(\Sigma)$ and $\pi(\mbox{int}(\Sigma))$. Thus, $\pi$ is 
an open map and $\pi(\mbox{int}(\Sigma))$ is a domain in $\Pi$. The compactness of $\Sigma$ implies that $\pi:\Sigma\rightarrow\overline{\Omega}$ is a covering map. Thus, we have

\begin{proposition}\label{p1}
Let $x:\Sigma\rightarrow\l^3$ be  a compact spacelike surface whose boundary 
$\Gamma$ is a graph over the boundary of a domain $\Omega\subset\r^2$. Then $x(\Sigma)$ is  a graph
over $\Omega$. 
\end{proposition}

We define the first and the second fundamental forms of $x$ as 
$${\rm I}=\sum_{i j}g_{ij} dx_i\ dx_j,\hspace*{1cm}{\rm II}=\sum_{i j} h_{ij} dx_i\ dx_j.$$
The mean curvature $H$ and the Gaussian curvature $K$ are given by 
$$2H=\frac{h_{22} g_{11}-2h_{12} g_{12}+h_{11}g_{22}}{\mbox{det}(g_{ij})},
\hspace*{1cm}K=-\frac{\mbox{det}(h_{ij})}{\mbox{det}(g_{ij})}.$$
If  $\Sigma$ is the graph of a function $x_3=u(x_1,x_2)$ defined over a domain $\Omega$,
 the spacelike condition implies$|Du|<1$ and
 the mean curvature $H$ is expressed by 
\begin{equation}\label{exp}
(1-|Du|^2)\sum_{i=1}^2 u_{ii}+\sum_{i,j=1}^2 u_i u_j u_{ij}= 2 H(1-|Du|^2)^{\frac{3}{2}}.
\end{equation}
This equation can alternatively be written in divergence form
$$\mbox{div}\left(\frac{Du}{\sqrt{1-|Du|^2}}\right)
=2 H.$$
Equation (\ref{exp}) is  of quasilinear elliptic type and the Hopf lemma can be applied. As a consequence, the solutions of the Dirichlet problem for the constant mean curvature equation are unique. See \cite{gt}.

 If  a  surface of revolution is parametrized by (\ref{form}), $f\in C^0([r,R])\cap C^2(r,R)$, 
the spacelike condition is equivalent to 
$$f'^2<1.$$
The computation of the mean curvature $H$ yields 
$$H=\frac{t f''(t)+(1-f'(t)^2)f'(t)}{2t(1-f'(t)^2)^{\frac32}}$$
with respect to the future-directed orientation. 
A first integral is obtained by 
$$\frac{d}{dt}\left(Ht^2-\frac{t f'(t)}{\sqrt{1-f'(t)^2}}\right)=0.$$
Thus, the quantity inside the parentheses is a constant $c$:
\begin{equation}\label{revo}
Ht^2-\frac{t f'(t)}{\sqrt{1-f'(t)^2}}=c.
\end{equation}
Equation (\ref{exp}) may be considered the Euler-Lagrange equation for critical points of the Minkowski area functional 
$\int_{\Omega}u\  dx_1\ dx_2=\mbox{ constant}$, with respect to strictly spacelike interior variations. We end 
this section showing that Equation (\ref{revo}) can directly derive   from volume-surface area considerations. Let  $\Sigma$ 
be a spacelike surface of revolution in $\l^3$ obtained by 
rotating the  curve $x_1=g(x_3)$ with respect to the $x_3$-axis. Assume that the profile curve
 has  fixed endpoints $r=g(a)$ and $R=g(b)$, $r<R$. The surface area  and the volume of $\Sigma$ are respectively 
$$A(\Sigma)=2\pi\int_r^R g(x_3)\sqrt{g'(x_3)^2-1}\ dx_3,\hspace*{1cm}V(\Sigma)=\pi\int_r^R g(x_3)^2\ dx_3.$$
We seek the surface which encloses a fixed volume $V(\Sigma)$ such that the surface area $A(\Sigma)$ is a critical 
point for any spacelike variation of $\Sigma$. 
Neglecting $\pi$ in the formula, we have to extremize the functional 
$$J=\int_r^R \left(2 g(x_3)\sqrt{g'(x_3)^2-1}-\lambda g(x_3)^2\right)\ dx_3=: \int_r^R F(g(x_3),g'(x_3))\ dx_3,$$
where $\lambda$ is the Lagrange multiplier.
We extremize this integral noting that the integrand is  independent of the variable $x_3$. The usual Euler-Lagrange argument says   that  there exists a constant $\kappa$ such that  the function $g$ satisfies  $F-g'\frac{\partial F}{\partial g'}=\kappa$. This gives
\begin{equation}\label{lagrange}\frac{2g}{\sqrt{g'^2-1}}-\lambda g^2=\kappa.\end{equation}
By considering $f=f(x_1)$ the inverse of the function $g$, Equation (\ref{lagrange}) writes as 
$$\frac{x_1 f'}{\sqrt{1-f'^2}}-\frac{\lambda}{2}x_1^2=\frac{\kappa}{2}$$
that coincides with (\ref{revo}) by taking $H=\lambda/2$ and $c=\kappa/2$.

%%%%%%%%%%%%%%%%%%%%%%%%%%%%%%%
\section{Proof of Theorem \ref{t1}}\label{s3}
%%%%%%%%%%%%%%%%%%%%%%%%%%%%%%%

Denote $f(t)=f(t;H,c)$ the solution of (\ref{revo}), emphasing its dependence on the values $H$ and $c$. 
It follows that $f(t;-H,-c)=-f(t;H,c)$.  Without loss of generality, we assume in this article that   $H\geq 0$. 
Moreover
\begin{equation}\label{ecuacion}
f'(t)=\frac{H t^2-c}{\sqrt{t^2+(H t^2-c)^2}}.
\end{equation}
If we denote $h(s)=f'(s)$, then
$$f(t)=a+\int_{r}^{t}h(s)\ {\rm ds}.$$
Since (\ref{revo}) is defined provided $t\not=0$ and the right-side hand in (\ref{ecuacion}) is 
a continuous function, the solutions of the differential equation (\ref{ecuacion}) are defined for any $t>0$. 
Moreover, the function $h(t)=h(t;H,c)$ is strictly non-increasing on $c$, and so, 
if $c_1<c_2$, then 
$$f(t;H,c_1)>f(t;H,c_2).$$
On the other hand,  if  $f(R)=b$, there exists $\xi\in [r,R]$ such that $h(\xi)(R-r)=b-a$. Then
$$\frac{|b-a|}{R-r}=|h(\xi)|<1,$$
and  this is a necessary condition for the existence of a spacelike surface of revolution spanning
$\Gamma(r,a)\cup\Gamma(R,b)$. This proves (i)$\Rightarrow$(ii) in Theorem \ref{t1}.
We prove the converse. Letting $c\rightarrow\pm\infty$, we obtain
$$\lim_{c\rightarrow+\infty}a+h(\xi)(R-r)=a-(R-r).$$
$$\lim_{c\rightarrow-\infty}a+h(\xi)(R-r)=a+(R-r).$$
By using the dependence of parameters for solutions for equation  (\ref{revo}), 
there exists a real number $c$ such that $f(R;H,c)=b$. This yields the desired solution.

We show that   at $t=0$ the surface presents a conical-type singularity unless that it is a planar plane or a hyperbolic plane. We have to prove that
$$\lim_{t\rightarrow 0}f'(t)^2=1.$$
When  $c=0$, it is possible to integrate (\ref{ecuacion}): if $H=0$,  the function  $f$ is a constant, that is, the surface is a horizontal plane; if $H\not=0$, then we obtain up constants that  $f(t)=\sqrt{1+H^2 t^2}/H$: this 
 surface  describes a hyperbolic plane and it is regular at $t=0$.

If $c\not=0$, then 
$$\lim_{t\rightarrow 0}f'(t)=-\frac{c}{|c|}.$$
Thus, 
\begin{enumerate}
\item If $c<0$, $\lim_{t\rightarrow 0}f'(t)=1$, and the surface is tangent to the upper light cone at $(0,f(0))$.
\item If $c>0$, $\lim_{t\rightarrow 0}f'(t)=-1$, and the surface is tangent to the lower light cone at $(0,f(0))$.
\end{enumerate}

Finally,  when $H\not=0$ the surface is asymptotic to a light cone at infinity provided
$$\lim_{t\rightarrow\infty}\frac{f(t)}{t}=\pm 1.$$
 L'Hôpital theorem yields  
$$\lim_{t\rightarrow\infty}\frac{f(t)}{t}=\lim_{t\rightarrow\infty}f'(t)=1.$$
Note that when $H=0$, $\lim_{t\rightarrow\infty}\frac{f(t)}{t}=0$. This completes the proof of Theorem \ref{t1}.

\begin{remark}\label{infinity}
{\rm The behaviour at infinity of a spacelike surface $u=u(x_1,x_2)$  
in Minkowski space $\l^3$ can 
be described by blowing $u$ down. See \cite{t}. We define the projective boundary value of $u$ at infinity by
$$V_u(x)=\lim_{s\rightarrow+\infty}\frac{u(s x)}{s},\hspace*{1cm}x=(x_1,x_2).$$
 Treibergs uses the concept of projective boundary value
to prove that the surface $u$ is asymptotically lightlike (see also \cite{ct}). For our  surfaces $f=f(t;H,c)$, we have 
$V_f=0$ if $H=0$, and $V_f(x)=|x|$ if $H\not=0$.}

\end{remark}

Combining  Theorem \ref{t1} and Proposition \ref{p1} we obtain immediately

\begin{corollary}\label{co} Let $\Gamma_1\cup \Gamma_2\subset\l^3$ be two concentric circles in parallel  planes. 
Let $\Sigma$ be a spacelike compact surface spanning $\Gamma_1\cup \Gamma_2$. If the   mean curvature of the 
surface is constant, then $\Sigma$ is a surface of revolution.
\end{corollary}

\begin{proof} We only point out that if $\Gamma$ is a closed planar curve included in  a spacelike surface, the plane containing 
$\Gamma$ must be spacelike. Thus and after a rigid motion of $\l^3$, we can suppose that 
$\Gamma_1\cup\Gamma_2$ lie in horizontal planes.
\end{proof}

\begin{remark}{\rm A theorem due to Shiffman states 
that a minimal surface in Euclidean space bounded by two circles in parallel planes must be foliated by circles 
in parallel planes \cite{sh}. In this sense, Corollary \ref{co} is a  partial  version of  Shiffman's theorem 
in the Lorentzian setting, with
the difference that in our case the circles are concentric and the mean curvature is a non-zero number. 
On the other hand, it has been proved that  a constant mean curvature spacelike surface 
foliated by circles  is a surface of revolution (if $H\not=0$) or it  is 
a Lorentzian catenoid or a Riemann-type surface (if $H=0$). See \cite{lo1,lo2}.}
\end{remark}

In the rest of this section,  we describe spacelike surfaces of revolution with constant mean curvature. We distinguish the cases
 $H=0$ and $H\not=0$.

%%%%%%%%%%%%%%%%%%%%%%%%%%%%%%%
\subsection{Case $H=0$}
%%%%%%%%%%%%%%%%%%%%%%%%%%%%%%%

\begin{figure}[h]
\includegraphics[width=9cm,height=5cm]{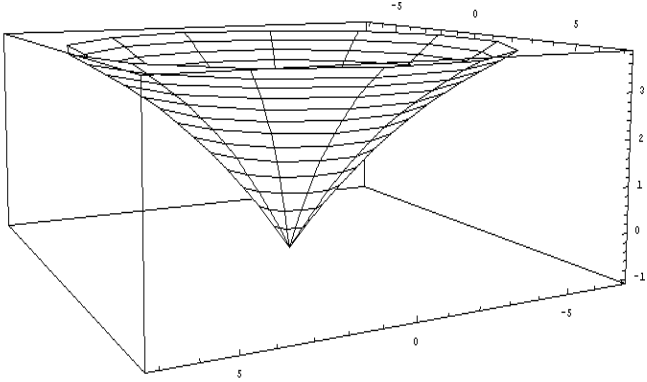}
\caption{A maximal surface: $f=f(t;0,3)$ with $f(1)=0$ and $0\leq t\leq 7$.}\label{figure1}
\end{figure}

It follows from Equation (\ref{revo}) that  the 
function $f$  satisfies
\begin{equation}\label{maximal}
\frac{t f'(t)}{\sqrt{1-f'(t)^2}}=c.
\end{equation}
For $c=0$, $f$ is a constant function and the surface is a horizontal planar domain. Assume $c\not=0$. A 
simple integration gives 
$$f(t;0,c)=  c \mbox{ arcsinh}{\left(\frac{t}{c}\right)}+d,\hspace*{1cm}d\in\r.$$
If we add the condition  $f(r)=a$, then 
$$f(t;0,c)= c\left(\mbox{arcsinh}\left(\frac{t}{c}\right)-\mbox{ arcsinh}\left(\frac{r}{c}\right)\right)+a.$$
See Figure \ref{figure1}.

\begin{remark} {\rm 
For a  maximal surface, the maximum principle  gets immediately that $\Sigma$ is included in 
the convex hull of its boundary. Thus, if the boundary of $\Sigma$ are two closed curves in parallel planes, 
the surface is included in the slab determined by both planes. When $\partial\Sigma$ are two concentric circles, the Alexandrov method of reflection by vertical planes
assures then $\Sigma$  is rotational symmetric. Moreover, Proposition \ref{p1} shows that the surface is a graph over 
a annular domain of $\r^2$. }
\end{remark}

%%%%%%%%%%%%%%%%%%%%%%%%%%%%%%%
\subsection{Case  $H\not=0$}
%%%%%%%%%%%%%%%%%%%%%%%%%%%%%%%
We analyse the three families of surfaces of revolution with constant mean curvature according 
to the sign of the parameter $c$.

%%%%%%%%%%%%%%%%%%
\noindent{\bf 1.  Case $c=0$}
%%%%%%%%%%%%%%%%%%

A direct integration gives
$$f(t)=\frac{\sqrt{1+H^2 t^2}}{H}+\mbox{constant},$$
and  the condition $f(r)=a$ gets 
$$f(t)=\frac{\sqrt{1+H^2 t^2}-\sqrt{1+H^2 r^2}}{H}+a.$$
This surface is a hyperbolic cap. Exactly, if $p=(0,0,\frac{\sqrt{1+H^2 r^2}}{H})$, the 
solution is a subset of the hyperbolic plane $\{x\in\l^3;\langle x-p,x-p\rangle=-\frac{1}{H^2}\}$.
 Therefore, 

\begin{proposition}\label{pr1} Let  $a, b, r $ and $R$ such that $0<b-a<R-r$. Then there exists a domain hyperbolic spanning $\Gamma(r,a)\cup \Gamma(R,b)$. The value $H$ of the mean curvature of this hyperbolic domain is 
$$H=H_0:=\frac{2(b-a)}{\sqrt{((R - r)^2-(b - a)^2) ((R + r)^2-(b - a)^2)}}.$$
\end{proposition}

%%%%%%%%%%%%%%%%%%
\noindent{\bf 2. Case $c<0$}
%%%%%%%%%%%%%%%%%%

For this case, the derivative $f'(t)$ does not vanish for any $t$ and $f'$ is positive. Thus
 $f$ is a strictly increasing function on $t$. As a consequence, any spacelike surface of revolution with constant mean curvature corresponding for $c<0$ and bounded by two concentric circles lies in the slab determined by the boundary planes.
 Denote $H_0$ the value obtained in Proposition \ref{pr1}.

\begin{proposition}\label{pr2} Let $a, b, r $ and $R$ such that $0<b-a<R-r$. Then for each positive real number $H$ such 
that 
\begin{equation}\label{condition2}
0<H<H_0\end{equation}
 there exists a unique spacelike surface of revolution $f=f(t;H,c)$,
with constant mean curvature $H$ spanning $\Gamma(r,a)\cup \Gamma(R,b)$ and where $c<0$ .
\end{proposition}

\begin{figure}
\includegraphics[width=9cm,height=5cm]{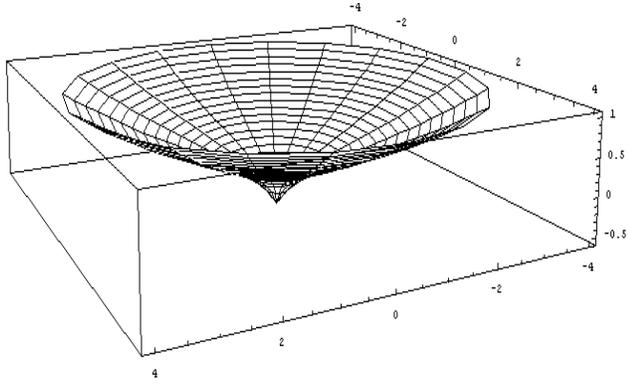}
\caption{The solution $f=f(t;\frac{1}{10},-\frac14)$ with $f(1)=0$ and $0\leq t \leq 4$.}\label{figure3}
\end{figure}

\begin{proof} Let $H$ be a number under the condition (\ref{condition2}). A straightforward computation  gives
$$a+\frac{\sqrt{1+H^2 R^2}-\sqrt{1+ H^2 r^2}}{H}<b.$$
We know that $f(R)=a+h(\xi)(R-r)$ for some $\xi$ with $r\leq\xi\leq R$.
The function $h(t)$ has a unique minimum at $t=\sqrt{-c/H}$. Thus
$$f(R;H,c)\geq a+ (R-r)h(\sqrt{-c/H})=a+(R-r)\frac{-2c}{\sqrt{4c^2-\frac{c}{H}}}.$$
As
$$\lim_{c\rightarrow -\infty}\frac{-2c(R-r)}{\sqrt{4c^2-\frac{c}{H}}}=R-r,$$
if $b<a+(R-r)$ and taking $c_0$ sufficiently close to  $-\infty$, we can find $c_0$ such that the corresponding solution 
for (\ref{ecuacion}) satisfies $f(R;H,c_0)>b$.

On the other hand, the theorem on continuity of parameters for  solutions of an
ordinary differential equation implies
$$f(R;H,c)\rightarrow f(R;H,0)=a+\frac{\sqrt{1+H^2 R^2}-\sqrt{1+ H^2 r^2}}{H},$$
as $c\rightarrow 0^-$. Continuity assures again the existence of a number $c_1$, $c_1<0$ such that 
$f(R;H,c_1)=b$. See Figure \ref{figure3}.
\end{proof}

%%%%%%%%%%%%%%%%
\noindent{\bf Case 3. $c>0$}
%%%%%%%%%%%%%%%%

When $c$ is a positive number, we have

\begin{proposition}\label{pr3} Let $a, b, r $ and $R$ such that $0\leq |b-a|<R-r$. Then for each positive real number $H$ such 
\begin{equation}\label{condition3}
H>H_0\end{equation}
 there exists a unique spacelike surface of revolution $f=f(t;H,c)$,
with constant mean curvature $H$ spanning $\Gamma(r,a)\cup \Gamma(R,b)$ and where $c>0$ .
\end{proposition}

\begin{proof} Let $f(t;H,c)$ be a solution of (\ref{revo}) with $f(a)=r$. Now condition (\ref{condition3}) implies that 
$$b<a+\frac{\sqrt{1+H^2 R^2}-\sqrt{1+ H^2 r^2}}{H}.$$
Because $c$ is a positive number, the tangent function $h$ 
is strictly increasing on $t$. Thus, for 
each $t\in [r,R]$, $h(r)\leq h(t)\leq h(R)$. It follows 
$$f(R;H,c)\leq a+h(R;H,c)(R-r)$$
and $\lim_{c\rightarrow \infty} a+h(R;H,c)(R-r)=a-(R-r)$. By continuity, there
exists $c_2>0$ and sufficiently big such that 
$f(R;H,c_2)<b$. On the other hand, 
and using  the theorem on continuity of parameters again, 
$\lim_{c\rightarrow 0}f(R;H,c)=f(R;H,0)=a+\frac{\sqrt{1+H^2 R^2}-\sqrt{1+ H^2 r^2}}{H}$. 
By continuity, there exists $c_3>0$ such that
$f(R;H,c_3)=b$.
\end{proof}

\begin{figure}
\includegraphics[width=9cm,height=5cm]{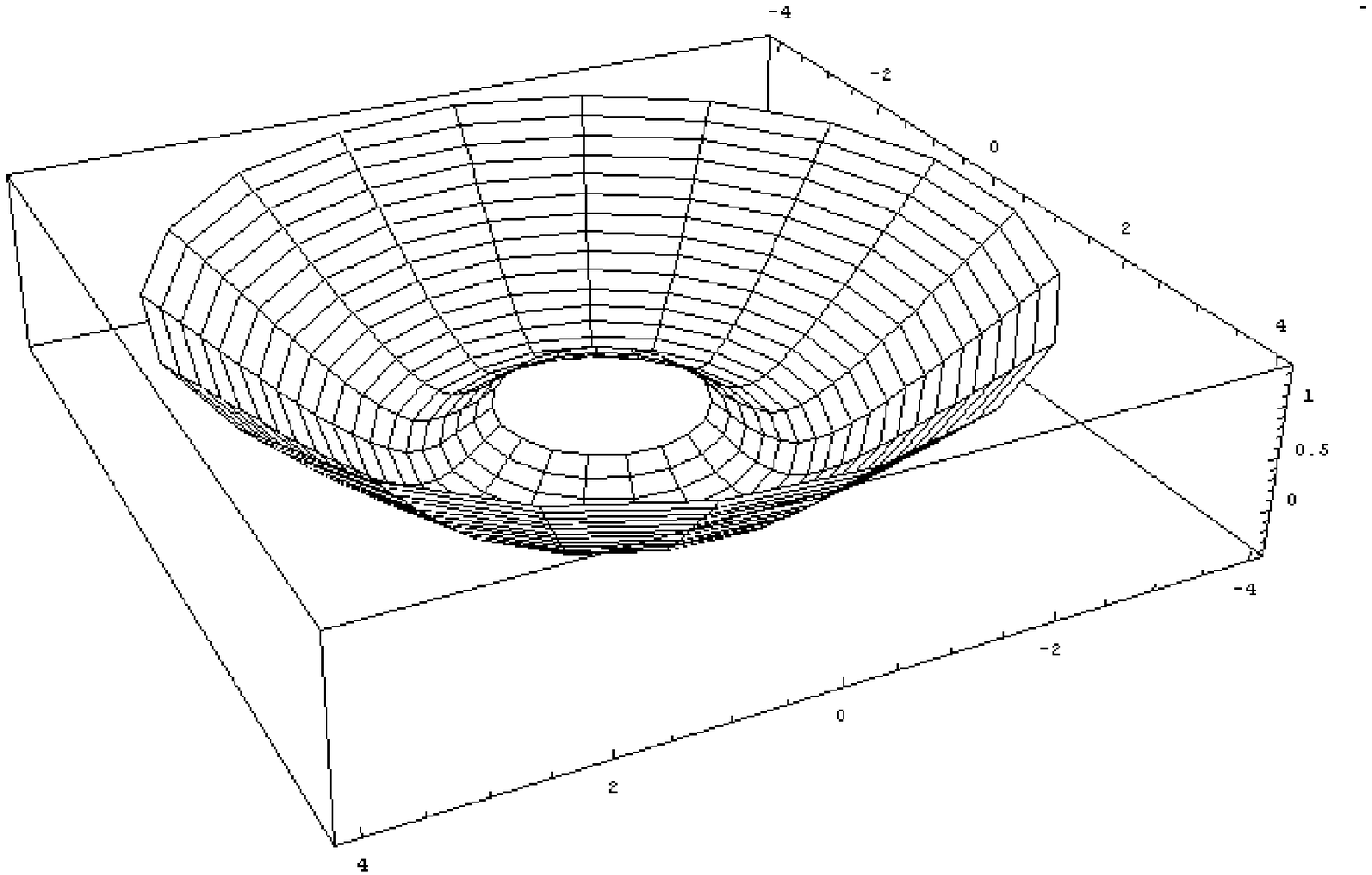}
\caption{The solution $f=f(t;1,3)$ with $f(1)=0$ and $1\leq t\leq 4$.}\label{figure41} 
\end{figure}

\begin{figure}
\includegraphics[width=9cm,height=5cm]{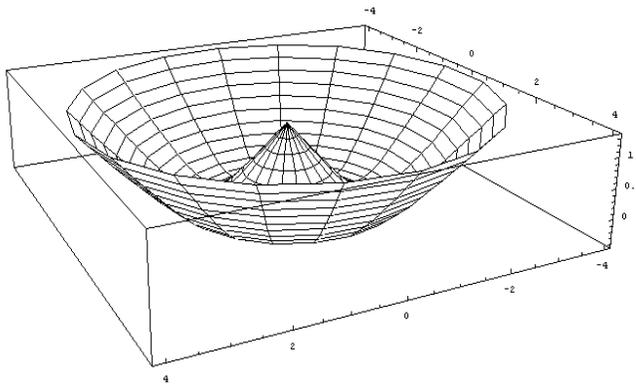}
\caption{The solution $f=f(t;1,3)$ with $f(1)=0$ and $0\leq t \leq 4$.}\label{figure42} 
\end{figure}

The function $h$  is increasing and vanishing at $t=\sqrt{c/H}$. Thus, 
$f$ is convex with a unique minimum. See Figure \ref{figure41}. 
Moreover, the case $a=b$ can be solved for some $c>0$. As a consequence, and taking values of $r<\sqrt{c/H}<R$, it is possible to find surfaces 
bounded by circles in parallel planes and not included in the slab determined by the boundary planes. See Figure \ref{figure42}.

%\begin{figure}
%\centereps{9cm}{5cm}{c:/articulos/annulus/lopez2.eps x=9cm y=5cm}
%\caption{The hyperbolic disk $f=f(t;2,0)$ with $0\leq t\leq 4$.}\label{figure2}
%\end{figure}

%%%%%%%%%%%%%%%%%%%%%%%%%%%%%%
\section{Proofs of Theorems \ref{t2} and \ref{t3}} \label{s5}
%%%%%%%%%%%%%%%%%%%%%%%%%
 
In this section we prove Theorems \ref{t2} and \ref{t3} using  the concept the flux of a curve in a surface and 
the tangency principle. The flux of a surface is used in a variety of problems in the theory of the constant mean curvature surfaces in Euclidean space
\cite{kk,kks}.

First, we introduce the concept of flux of a surface. 
Let $\e_3=(0,0,1)$ . Consider  $x:\Sigma\rightarrow \l^3$  a spacelike immersion with constant mean curvature $H$. 
Define the 1-form
$$\omega_p(v)=H \langle x\wedge dx_p(v),\e_3\rangle +\langle N(p)\wedge d x_p(v),\e_3\rangle,\hspace*{1cm}p\in \Sigma$$
where $v$ is a vector tangent to $\Sigma$ at $p$ and $\wedge$ is the cross-product in $\l^3$. The constancy of the 
the mean curvature implies that  $\omega$ is a closed form. It follows from Stokes 
theorem that given a 1-cycle  $\Gamma$  on 
 $\Sigma$ bounding  an open $Q\subset\Sigma$,  the expression
\begin{equation}\label{flux}
{\rm Flux}(\Gamma)=H\int_{\Gamma}\langle x\wedge\tau,\e_3\rangle{\rm ds}+\int_{\Gamma}\langle \nu,\e_3\rangle{\rm ds}
\end{equation}
depends only on the homology class of $\Gamma$ on $\Sigma$. Here $\tau$ is a unit tangent vector to $\partial\Sigma$, 
and $\nu$ is the unit conormal vector of  $Q$ along $\Gamma$ such that $N\wedge \tau=\nu$. Note that 
the first summand in (\ref{flux}) is $2H$ times the algebraic area of the orthogonal projection of 
$\Gamma$ on the plane $x_3=0$. The number ${\rm Flux}(\Gamma)$ is called 
the flux of $\Gamma$. For example, the flux of a null-homologous cycle is zero. Formula (\ref{flux}) 
can be viewed as a measure of the forces of the surface tension of $ \Sigma$ 
that act along its boundary and  the pressure forces 
that act on  $\partial \Sigma$.

Consider now the surfaces of revolution obtained in Section \ref{s3}.
Denote $f=f(t;H,c)$ the solution of (\ref{revo}), $r\leq t\leq R$, with  boundary condition 
$f(r)=a$,  and let 
$\Gamma(r)=\{(x_1,x_2,f(r));x_1^2+x_2^2=r^2\}$. 
Here  
$$\int_{\Gamma(r)}\langle\nu,\e_3\rangle{\rm \ ds}=\int_{\Gamma(r)}\frac{-f'(r)}{\sqrt{1-f'(r)^2}}{\rm \ ds}=
-2\pi (H r^2-c).$$
Thus
$${\rm Flux}(\Gamma(r))=2\pi c.$$
As conclusion, if we fix real numbers $H,r,a $ and $\lambda$,
{\it there is a  rotational symmetric spacelike graph on  $\{(x_1,x_2);x_1^2+x_2^2>r^2\}$ 
spanning $\Gamma(r,a)$, with  constant mean curvature $H$ and such that ${\rm Flux}(\Gamma(r,a))=\lambda$}.

Finally,  we will state the tangency principle for spacelike surfaces with constant mean curvature.
Let   $u$ and $v$ be two functions that are local expressions of two spacelike surfaces $\Sigma_u$ and $\Sigma_v$ 
of  $\l^{3}$. If $\Sigma_u$ and $\Sigma_v$ have a common point $p=(p_1,p_2,p_3)$ where they are tangent, 
we will say that $\Sigma_u$ lies above $\Sigma_v$ near $p$ when $u\geq v$ on a 
certain neighborhood of the point $(p_1,p_2)$. Let us assume  that $\Sigma_u$ and 
$\Sigma_v$ have the same constant mean curvature $H$. Since Equation (\ref{exp}) is  of quasilinear elliptic type, the difference function $u-v$ satisfies a {\it linear} elliptic 
equation on a neighborhood of $(p_1,p_2)$ and the Hopf 
maximum principle for linear elliptic equations  can be applied to $u-v$ 
(see \cite[Th. 9.2]{gt}). The same holds if $p$ is a common boundary point with the extra hypothesis that $\partial\Sigma_1$ and  $\partial\Sigma_2$ are 
tangent at $p$.        
  
Consequently, we have proved the following result.

\begin{proposition}[Tangency principle] Let $\Sigma_1$ and $\Sigma_2$ be two spacelike surfaces 
(possibly with boundary) in $\l^{3}$ with the same constant mean curvature with respect to the future-directed orientation. Suppose  they 
intersect tangentially at a point $p$. If $\Sigma_1$ is above 
$\Sigma_2$, then $\Sigma_1=\Sigma_2$ locally around $p$ in a neighborhood of $p$ if one of the
following hypotheses holds:
\begin{enumerate}
\item $p$ is an interior point of $\Sigma_1$ and $\Sigma_2$.
\item $p$ is a boundary point of $\Sigma_1$ and $\Sigma_2$ and 
$\partial\Sigma_1$ and $\partial\Sigma_2$ are tangent at $p$.
\end{enumerate}
 \end{proposition}

In either case, by analyticity of solutions of elliptic equations, we conclude that 
$\Sigma_1$ and $\Sigma_2$ coincide whenever they are simultaneously defined. 

Now we are in position to prove Theorems \ref{t2} and \ref{t3}.

\begin{proof}[of Theorem \ref{t2}] Let  $u\in 
C^2(\Omega)\cap C^0(\Omega)$  be a  solution of the constant mean curvature equation (\ref{exp}) in $\Omega$
 with  $u=a$ on 
$\partial\Omega$. Let $\lambda=Flux(\Gamma(r,a))$ be the flux of $\Gamma(r,a)$ in $\Sigma_1$, the graph of $u$.
We consider $f=f(t;H,c)$, with $t\in [r,\infty)$,  $f(r)=a$ and $c$ satisfying $2\pi c=\lambda$.  
Denote $\Sigma_2$ the graph  of the function $f$. 

In particular, the flux of 
the boundary curve $\Gamma(r,a)$ agrees in both surfaces. Since the first summand in (\ref{flux})  is the same 
for both surfaces, we conclude that 
$$\int_{\Gamma(r,a)}\langle \nu_1,\e_3\rangle {\rm \ ds}=\int_{\Gamma(r,a)}\langle \nu_2,\e_3\rangle {\rm \ ds},$$
where $\nu_i$ means the conormal on the surface  $\Sigma_i$, $i=1,2$.

We prove that  $\Sigma_1=\Sigma_2$. We move upwards the surface $\Sigma_2$ and 
put $\Sigma_1(t)=\Sigma_1+t\e_3$. Recall that both surfaces are asymptotically lightlike (Remark 
\ref{infinity}). Therefore, for $t>0$ big enough, $\Sigma_1(t)$ 
does not intersect $\Sigma_2$. Let us descend $\Sigma_1(t)$ until the first time $t_0$ that touches $\Sigma_2$, that is, 
$t_0=\inf\{t>0;\Sigma_1(t)\cap \Sigma=\emptyset\}$. We want to prove that $t_0=0$. If $t_0>0$ we have two possibilities.
First, if $\Sigma_1(t_0)\cap\Sigma_2\not=\emptyset$, then both surfaces have a contact {\it interior point}. Because both 
surfaces have the same mean curvature, the tangency principle implies 
that $\Sigma_1(t_0)=\Sigma_2$, which is impossible. If $\Sigma_1(t_0)\cap\Sigma_2=\emptyset$, then 
$\Sigma_1(t_0)$ and $\Sigma_2$ have a contact at infinity. Thus it is possible to choose 
$\epsilon$ such that  $\Gamma=\Sigma_1(t_0-\epsilon)\cap \Sigma_2$ is a 1-dimensional analytic 
curve. By the maximum principle, there is a unique graph with constant mean curvature spanning a given 
closed curve and thus, $\Gamma$ is a analytic Jordan curve  in the same homology class 
of $\Gamma(r,a)$.  Along 
$\Gamma(r,a)$, we have $\langle\nu_2,\e_3\rangle<\langle\nu_1,\e_3\rangle$, and so, 
$$\int_{\Gamma(r,a)} \langle\nu_2,\e_3\rangle{\rm \ ds}<\int_{\Gamma(r,a)} \langle\nu_1,\e_3\rangle{\rm \ ds}.$$
However this contradicts the fact the fluxes of $\Gamma(r,a)$ coincide for $\Sigma_1(t_0-\epsilon)$ and 
$\Sigma_2$.

This contradiction implies that the assumption of $t_0>0$ is false, and we have $t_0=0$. By repeating the argument, either $\Sigma_1$ and $\Sigma_2$ have 
a contact {\it boundary point}, and which case the tangency principle gets $\Sigma_1=\Sigma_2$, or 
$\langle\nu_2,\e_3\rangle<\langle\nu_1,\e_3\rangle$ along $\Gamma(r,a)$. Since 
the fluxes of $\Gamma(r,a)$ for $\Sigma_1$ and 
$\Sigma_2$ are the same along this curve, we obtain an absurd again. This completes the proof of Theorem \ref{t2}.
\end{proof}

\begin{proof}[of Theorem \ref{t3}]

Let  $\Sigma$ be an entire surface with constant mean curvature $H$ defined by the graph of a  function
$u$ and without loss of generality, we assume  $u$ is defined  on $\r^2\setminus\{(0,0)\}$ where the 
origin the is the singularity of $\Sigma$. 
When $H=0$, Theorem \ref{t3} was proved in \cite{ec}. Thus, assume $H\not=0$. After a horizontal symmetry, we can assume that $H$ is positive.

 The method of proof is similar as in  Theorem \ref{t2}. Fix $r>0$ and let us 
consider $f=f(t;H,c)$, $t\geq r$, the rotational symmetric spacelike surface with constant mean curvature $H$ such that 
$f(r)=a$ and $c$ satisfying $2\pi c=\lambda$.  We extend $f$ to $(0,\infty)$. 
Denote $\Sigma_2$ the graph  of the function $f$ and $\Sigma_1=\Sigma$.
Proceeding similarly as in Theorem \ref{t2}, we 
conclude that $t_0=0$, that is, $\Sigma_2$ is over $\Sigma_1$ and they contact in the 
singularity. By reversing argument, we prove that $\Sigma_1$ is over $\Sigma_2$, and so, 
$\Sigma_1=\Sigma_2$.
\end{proof}

%%%%%%%%%%%%%%%%%%%%

\end{document}